\documentclass[a4paper,12pt]{article}
\title{Arcsine Law as the Classical Limit \\for interacting Fock spaces}
\author{Hayato Saigo\footnote{E-mail: h\_saigoh@nagahama-i-bio.ac.jp } 
\\ Nagahama Institute of Bio-Science and Technology \\  Nagahama 526-0829, Japan
}
\date{}

\usepackage{amsmath}
\usepackage{amsthm}
\usepackage{amssymb}
\usepackage{amsbsy}
\usepackage{amsfonts}
\usepackage{amstext}
\usepackage{amscd}
\usepackage[dvips]{epsfig}

\numberwithin{equation}{section}
\theoremstyle{plain}
\newtheorem{thm}{Theorem}[section]
\newtheorem{prop}[thm]{Proposition}

\newtheorem{lem}[thm]{Lemma}
\theoremstyle{definition}

\newtheorem{rem}[thm]{Remark}
\newtheorem{df}[thm]{Definition}
\newtheorem{nota}[thm]{Notation}

\begin{document}
\maketitle

\begin{abstract}
 In the present paper we discuss how to generalize ``Quantum-Classical Correspondence'' 
by means of the notion of interacting Fock spaces, 
which associates algebraic probability theory and the theory of orthogonal polynomials of probability measures. 
As an application we show that the Arcsine Law is ``Classical Limit'' for interacting Fock spaces corresponnding to certain kind of symmetric probability measures such as q-Gaussians. 
We also discuss the case of the exponential distribution as a simple example of asymmetric probability measures. 
\end{abstract}

\section{Introduction}
The distribution $\mu_{As}$ defined as
\[
\mu_{As} (dx)=\frac{1}{\pi} \frac{dx}{\sqrt{2-x^2}}\:\:\: (-\sqrt{2}<x<\sqrt{2}).
\]
is called the (normalized) \textbf{Arcsine Law}, which plays lots of crucial roles both in pure and applied probability theory. The $n$-th moment $M_n:=\int_{\mathbb{R}}x^n \mu_{As}(dx)$ is given by
\[
M_{2m+1}=0,\:\:\:M_{2m}=\frac{1}{2^m}\binom{2m}{m}.
\]
The moment problem for the Arcsine law is determinate, that is, the moment sequence $\{M_n\}$ characterizes $\mu_{As}$.
In \cite{SAI} we have proved 
that the Arcsine Law appears as the \textbf{Classical Limit} of quantum harmonic oscillator, in the framework of 
algebraic probability thoery (also known as ``noncommutative probability theory'' or ``quantum probability theory''). 

The purpose of this paper is to extend this ``quantum-classical correspondence'' in general interacting Fock spaces \cite{A-B}. 
It implies asymptotic behavior of orthogonal polynomials for certain kind of symmetric probability measures as we see in Section 4. 
In section 5 we also discuss the ``classical limit'' for the case of exponential distribution, as a simple example of asymmetric probability measures.

\section{Basic Concepts}

Let $\mathcal{A}$ be a $\ast$-algebra. We call a linear map $\varphi :\mathcal{A}\rightarrow \mathbb{C}$ a state on $\mathcal{A}$ if it satisfies
\[
\varphi(1)=1,\:\:\:\varphi (a^{\ast}a)\geq 0.
\]
A pair $(\mathcal{A}, \varphi)$ of a $\ast$-algebra and a state on it is called an algebraic probability space.
Here we adopt a notation for a state $\varphi :\mathcal{A}\rightarrow \mathbb{C}$, an element $X\in \mathcal{A}$ and a probability distribution $\mu$ on $\mathbb{R}$.
\begin{nota} We use the notation 
$X\sim_{\varphi} \mu $
when $\varphi(X^m)=\int_{\mathbb{R}}x^m\mu (dx)\:\:\: $for all $m \in \mathbb{N}$.
\end{nota}
\begin{rem}
Existence of $\mu$ for $X$ which satisfies $X\sim_{\varphi} \mu$ always holds. 
\end{rem}

\begin{df}[Jacobi sequence] A sequence $\{\omega_n\}$ is called a Jacobi sequence if it satisfies one of the conditions below:
\begin{itemize}
\item (finite type) There exist a number $m$ such that $\omega_n>0$ for $n< m$ and $\omega_n=0$ for $n\geq m$;
 \item (infinite type) $\omega_n>0$ for all $n$.
 \end{itemize}
\end{df}

\begin{df}[Interacting Fock space] Let $\{\omega_n\}$ be a Jacobi sequence. An interacting Fock space $\Gamma_{\{\omega_n\}}$ is a triple 
$(\Gamma(\mathbb{C}), a, a^{\ast})$ where $\Gamma(\mathbb{C})$ is a Hilbert space 
$\Gamma(\mathbb{C}):=\oplus^{\infty}_{n=0} \mathbb{C}\Phi _n$ with inner product given by $\langle\Phi_n, \Phi_m\rangle=\delta _{n,m}$, 
and $a, a^{\ast}$ are operators defined as follows:
\[
a\Phi_0 =0, \:\:\:a\Phi_n=\sqrt{\omega_n}\Phi_{n-1} (n\geq 1)
\]
\[
a^{\ast}\Phi _n=\sqrt{\omega_{n+1}}\Phi_{n+1}. \]
\end{df}
Let $\mathcal{A}$ be the ${\ast}$-algebra generated by $a$, 
and $\varphi_n$ be the state defined as $\varphi_n(\cdot):=\langle\Phi_n, (\cdot)\Phi_n\rangle$. 
Then $(\mathcal{A}, \varphi_n)$ is an algebraic 
probability space. 

The interacting Fock space corresponding to $\omega_n=n$ is called ``Quantum Harmonic Osillator''. For quantum harmonic oscillator, it is well known that 
\[
a+a^{\ast}
\]
represents the ``position'' and that 
\[
a+a^{\ast}\sim _{\varphi_0} \frac{1}{\sqrt{2\pi}}e^{-\frac{1}{2}x^{2}}dx.
\]
That is, in $n=0$ case, the distribution of position is Gaussian. 

On the other hand, the asymptotic behavior of the distributions of position as $n$ tends to infinity is nontrivial. In other words, what is the ``Classical limit'' of quantum harmonic oscillator?

\section{Quantum-Classical Correspondence}
This question, 
which is related to fundamental problems in Quantum theory and asymptotic analysis \cite{EZA}, 
was analyzed in \cite{SAI} from the viewpoint of noncommutative algebraic probability 
with quite a simple combinatorial argument. The answer for this question is that the ``Classical Limit'' for quantum harmonic oscillator is nothing but the Arcsine law. Here we generalize this result:

\begin{thm}\label{IFS} Let $\Gamma_{\{\omega_n\}}:=(\Gamma(\mathbb{C}), a, a^{\ast})$ be an interacting Fock space satisfying the condition
\[
\lim_{n\rightarrow \infty}\frac{\omega_{n+1}}{\omega_n}=1 
\]
and $\mu_N$ be a probability distribution on $\mathbb{R}$ such that 
\[
\frac{a+a^{\ast}}{\sqrt{2\omega_N}}\sim _{\varphi_N} \mu_N .
\]
Then $\mu_N$ weakly converges to $\mu_{As}$.
\begin{proof}
Since it is known that moment convergence implies weak convergence when the moment problem for the limit distribution is determinate (Theorem 2.5.5 in \cite{CHU}), 
we are going to show the limit moment is the moment of the Arcsine law.


First, it is clear that
\[
\varphi_N((\frac{a+a^{\ast}}{\sqrt{2\omega_N}})^{2m+1})=\langle\Phi_N,(\frac{a+a^{\ast}}{\sqrt{2\omega_N}})^{2m+1}\Phi_N \rangle=0
\]
since $\langle\Phi_N, \Phi_M\rangle=0$ when $N\neq M$. 

To consider the moments of even degrees, we introduce the following notations:
\begin{itemize}
 \item $\Lambda^{2m}:=\{\text{maps from $\{1,2,..., 2m\}$ to $\{1,\ast\}$}\}$, 
 \item $\Lambda^{2m}_{m}:=\{\lambda \in \Lambda^{2m}; |\lambda^{-1}(1)|=|\lambda^{-1}(\ast)|=m\}$.
\end{itemize}

Note that the cardinality $|\Lambda^{2m}_{m}|$ equals to $\binom{2m}{m}$ because the choice of $\lambda$ is equivalent to the choice of $m$ elements which consist the subset $\lambda^{-1}(1)$ from 
$2m$ elements in $\{1,2,..., 2m\}$.

It is clear that for any $\lambda \notin  \Lambda^{2m}_{m}$ 
\[
\langle\Phi_N, a^{\lambda_1}a^{\lambda_2}\cdots a^{\lambda_{2m}}\Phi_N\rangle=0
\]
since $\langle\Phi_N, \Phi_M\rangle=0$ when $N\neq M$.

On the other hand, for any $\lambda \in \Lambda^{2m}_{m}$ 
\[
\frac{1}{\omega_N^m} \langle\Phi_N, a^{\lambda_1}a^{\lambda_2}\cdots a^{\lambda_{2m}}\Phi_N\rangle \:\rightarrow 1 \:\:\:\:(N\rightarrow \infty)
\]
holds since $\langle\Phi_N, a^{\lambda_1}a^{\lambda_2}\cdots a^{\lambda_{2m}}\Phi_N\rangle$ becomes the product of $2m$ terms having the form
$\sqrt{\omega_{N+k}}$ ($k$ is an integer and $-m+1\leq k\leq m$) and 
\[
\frac{\omega_{N+k}}{\omega_N} \rightarrow 1 \:\:\:\:(N\rightarrow \infty)
\]
by the assumption. Hence,
\[
 \varphi_N((\frac{a+a^{\ast}}{\sqrt{2\omega_N}})^{2m})=\langle\Phi_N,(\frac{a+a^{\ast}}{\sqrt{2\omega_N}})^{2m}\Phi_N \rangle 
\]
\begin{eqnarray*}
=&\frac{1}{2^m}\sum_{\lambda\in 
\Lambda^{2m}} \frac{1}{\omega_N^{m}}\langle\Phi_N, a^{\lambda_1}a^{\lambda_2}\cdots a^{\lambda_{2m}}\Phi_N\rangle & \\
\\
=&\frac{1}{2^m}\sum_{\lambda\in \Lambda^{2m}_{m}} \frac{1}{\omega_N^{m}}\langle\Phi_N, a^{\lambda_1}a^{\lambda_2}\cdots a^{\lambda_{2m}}\Phi_N\rangle  & 
\end{eqnarray*}
\[
\rightarrow \frac{1}{2^m}|\Lambda^{2m}_{m}|=\frac{1}{2^m}\binom{2m}{m} \:\:\:\:(N\rightarrow \infty).\\
\]
\end{proof}

\end{thm}

\section{Asymptotic behavior of Orthogonal Polynomials}

The theorem above has an interpretation in terms of orthogonal polynomials. To see this we review the relation between interacting Fock spaces, probability measures and orthogonal polynomials.

Let $\mu$ be a probability measure on $\mathbb{R}$ having finite moments. 
(For the rest of the present paper, we always assume that all the moments are finite.)
Then the space of polynomial functions is contained in the Hilbert space $L^2(\mathbb{R},\mu )$. 
A Gram-Schmidt procedure provides orthogonal polynomials which only depend on the moment sequence.

Let $\{p_n(x) \}_{n=0, 1, \cdots}$ be the monic orthogonal polynomials of $\mu$ such that 
the degree of $p_n$ equals to $n$. 
Then there exist sequences $\{\alpha_n\}_{n=0, 1, \cdots}$ 
and Jacobi sequence $\{\omega_n\}_{n= 1, 2, \cdots}$ such that 
\[
x p_n(x) 
= p_{n+1}(x) + \alpha_{n+1} p_n(x) + \omega_n p_{n-1}(x) \quad (p_{-1}(x)\equiv 0).
\]
$\alpha_{n}\equiv 0$ if $\mu$ is symmetric, i.e., $\mu(-dx)=\mu(dx)$. 

It is known that there exist an isometry $U:\Gamma_{\{\omega_n\}}\rightarrow L^2(\mathbb{R},\mu )$ through which we obtain 
\[
a+a^{\ast}+a^o \sim _{\varphi_N} |P_N(x)|^2\mu(dx)
\]
where $a^o$ is an operator defined by $a^o\Phi_n:=\alpha_{n+1}\Phi_n $
and $P_n$ denotes the normalized orthogonal polynomial of degree $n$ \cite{H-O}. 
Then Theorem \ref{IFS} implies the following: 
\begin{thm}
Let  $\mu$ be a symmmetric measure such that the corresponting Jacobi sequence $\{\omega_n\}$ satisfies 
\[
\lim_{n\rightarrow \infty}\frac{\omega_{n+1}}{\omega_n}=1 
\]
Then the measure $\mu_n$ defined as 
$\mu_n(dx):=|P_n(\sqrt{2\omega_n}x)|^2\mu(\sqrt{2\omega_n}dx)$ weakly converge to $\mu_{As}$.
\end{thm}
Since ``q-Gaussians''  ($-1< q\leq 1$, $q=1$ is Gaussian and $q=0$ is Wigner Semicircle Law), corresponding to $\omega_n=[n]_q:=1+q+q^2+\cdots +q^{n-1}$, satisfy the condition above, 
$\mu_{As}$ is turned out to be the Classical Limit of these measures.

In the next section we discuss the Classical Limit for the case of exponential distribution as an example of asymmetric measure.

\section{Exponential-Laguerre case}
Let $\mu$ be the exponential distribution, i.e., $\mu (dx):=e^{-x}dx$ $(x>0)$.Then
\[
x l_n(x) 
= l_{n+1}(x) + (2n+1)l_n(x) + n^2 l_{n-1}(x) \quad (l_{-1}(x)\equiv 0),
\]
holds, where $l_n$ denotes the Laguerre polynomial of $n$-th degree modified to be monic. Let us consider the interacting Fock space $\Gamma_{\{\omega_n\}}$ for $\omega_n=n^2$. 
As we have discussed, 
\[
a+a^{\ast}+a^o \sim _{\varphi_N} |L_N(x)|^2e^{-x}dx \:\:\:(x>0).
\]
where $L_n$ denotes the usual (normalized) Laguerre polinomial of $n$-th order.

Then we can calculate the ``Limit moment'' of $\mu_n(dx):=|L_n(nx)|^2ne^{-nx}dx$ $(x>0)$ in the spirit of the proof of Theorem \ref{IFS}. The result is:

\begin{prop}
\[
\lim_{N\rightarrow \infty}\varphi_N((\frac{a+a^{\ast}+a^{o}}{N})^m)=\sum_l 2^{m-2l}\binom{m}{m-2l}\binom{2l}{l}.
\]
\end{prop}

The right hand side of the proposition above is simplified as follows.

\begin{lem}
\[\sum_l 2^{m-2l}\binom{m}{m-2l}\binom{2l}{l}=\binom{2m}{m}.\]
\begin{proof}
Consider two sets of maps
\[L:=\{f:\mathbf{m}\rightarrow \mathbf{4} ; |f^{-1}(0)|=|f^{-1}(1)|\}\]
\[R:=\{\tilde{f}:\mathbf{2\times m}\rightarrow \mathbf{2} ; |\tilde{f}^{-1}(0)|=|\tilde{f}^{-1}(1)|\},\]
where $\mathbf{m}:=\{0,1,2,\dots,m-1\}$. Since we can construct an isomorphism between $L$ and $R$,  $|L|=|R|$. This is what to be proved. (This proof is obtained in discussion with Hiroki Sako).
\end{proof}
\end{lem}

It is easy to show that
\[
\binom{2m}{m}=\int_{0}^{4}x^m \frac{1}{\pi}\frac{dx}{\sqrt{4-(x-2)^2}},
\]
and hence we obtain the following theorem.

\begin{thm}
Let $L_n$ be the normalized Laguerre polynomial of $n$-th degree. Then 
$\mu_n(dx):=|L_n(nx)|^2ne^{-nx}dx$ $(x>0)$ weakly converge to 
\[
\frac{1}{\pi}\frac{dx}{\sqrt{4-(x-2)^2}}\:\:\:(0<x<4).
\]
\end{thm}
That is, the \textbf{Classical Limit of ``Laguerre oscillator'' is also the Arcsine Law} (just translated and dilated). 

\section*{Acknowledgments} 
The author is grateful to Prof. Marek Bo\.{z}ejko for advises and encouragements. He would like to thank Prof. Izumi Ojima and Mr. Kazuya Okamura  
for discussions on Quantum-Classical Correspondence. He deeply appreciates Hiroki Sako for collaboration.

\end{document}